\documentclass{article}

\usepackage{
 amsmath,
 amsxtra,
 amsthm,
 amssymb,
 etex,
 mathrsfs,
 mathtools,
 tikz-cd,
 xr,
 fullpage,
 enumerate,
 comment
 }
\usepackage[all]{xy}
\usepackage{hyperref}

\newtheorem{theorem}{Theorem}[subsection]
\newtheorem{lemma}[theorem]{Lemma}

\newtheorem{proposition}[theorem]{Proposition}
\newtheorem{corollary}[theorem]{Corollary}

\theoremstyle{remark}
\newtheorem{remark}[theorem]{Remark}

\newtheorem{defn}[theorem]{Definition}

\newcommand{\bd}{\begin{defn}}
\newcommand{\ed}{\end{defn}}
\newcommand{\bl}{\begin{lemma}}
\newcommand{\el}{\end{lemma}}
\newcommand{\bp}{\begin{proposition}}
\newcommand{\ep}{\end{proposition}}
\newcommand{\bt}{\begin{theorem}}
\newcommand{\et}{\end{theorem}}
\newcommand{\bc}{\begin{corollary}}
\newcommand{\ec}{\end{corollary}}
\newcommand{\br}{\begin{remark}}
\newcommand{\er}{\end{remark}}
\newcommand{\ba}{\begin{array}}
\newcommand{\ea}{\end{array}}
\newcommand{\bpf}{\begin{proof}}
\newcommand{\epf}{\end{proof}}

\newcommand{\Q}{\mathbb{Q}}
\newcommand{\Qp}{\mathbb{Q}_p}
\newcommand{\Z}{\mathbb{Z}}
\newcommand{\Zp}{\mathbb{Z}_p}
\newcommand{\Op}{\mathcal{O}}
\newcommand{\Ga}{\Gamma}
\newcommand{\ga}{\gamma}
\newcommand{\La}{\Lambda}
\newcommand{\la}{\lambda}

\DeclareMathOperator{\Gal}{Gal}
\DeclareMathOperator{\Hom}{Hom}

\DeclareMathOperator{\rank}{rank}
\DeclareMathOperator{\Ext}{Ext}
\newcommand{\Iw}{\mathrm{Iw}}

\newcommand{\lra}{\longrightarrow}
\newcommand{\ord}{\mathrm{ord}}
\newcommand{\ot}{\otimes}

\newcommand{\cyc}{\mathrm{cyc}}
\newcommand{\ps}[1]{[[ #1 ]]}

\newcommand{\ilim}{\displaystyle \mathop{\varinjlim}\limits}
\newcommand{\plim}{\displaystyle \mathop{\varprojlim}\limits}

\numberwithin{equation}{section}

\begin{document}
\title{On the growth of even $K$-groups of rings of integers in $p$-adic Lie extensions}
 \author{
  Meng Fai Lim\footnote{School of Mathematics and Statistics $\&$ Hubei Key Laboratory of Mathematical Sciences,
Central China Normal University, Wuhan, 430079, P.R.China.
 E-mail: \texttt{limmf@mail.ccnu.edu.cn}} }
\date{}
\maketitle

\begin{abstract} \footnotesize
\noindent Let $p$ be an odd prime number. In this paper, we study the growth of the Sylow $p$-subgroups of the even $K$-groups of rings of integers in a $p$-adic Lie extension. Our results generalize previous results of Coates and Ji-Qin, where they considered the situation of a cyclotomic $\Zp$-extension. Our method of proof differs from these previous work. Their proof relies on an explicit description of certain Galois group via Kummer theory afforded by the context of a cyclotomic $\Zp$-extension, whereas our approach is via considering the Iwasawa cohomology groups with coefficients in $\Zp(i)$ for $i\geq 2$. We should mention that this latter approach is possible thanks to the Quillen-Lichtenbaum Conjecture which is now known to be valid by the works of Rost-Voevodsky. We also note that the approach allows us to work with more general $p$-adic Lie extensions that do not necessarily contain the cyclotomic $\Zp$-extension, where the Kummer theoretical approach does not apply. Along the way, we study the torsionness of the second Iwasawa cohomology groups with coefficients in $\Zp(i)$ for $i\geq 2$. Finally, we give examples of $p$-adic Lie extensions, where the second Iwasawa cohomology groups can have nontrivial $\mu$-invariants.


\end{abstract}

\section{Introduction}

Iwasawa \cite{Iw59, Iw73} first proved a remarkable regularity on the variation of the Sylow $p$-subgroup of the class groups in the intermediate subfields of a $\Zp$-extension of a number field $F$. From an algebraic $K$-theoretical point of view, the class group is none other than the torsion subgroup of $K_0(\Op_F)$, where $\Op_F$ is the ring of integers of the number field $F$. It therefore seems natural to consider the analogous situation for the higher $K$-groups, considering that the higher $K$-groups $K_{m}(\Op_F)$ are known to encode important arithmetic information on the field $F$ (for instance, see \cite{Bo, C72, Gar, Kol, Sch, Sou, Wei05}). In this direction, Coates \cite{C72} established an analogue of Iwasawa's result for the orders of Sylow $p$-subgroups of the $K_2$-groups of rings of integers in a cyclotomic $\Zp$-extension, and this was subsequently extended to higher even $K$-groups by Ji-Qin \cite{JQ} (also over a cyclotomic $\Zp$-extension).

At the turn of the millennium, Coates and many others initiated a major study on towers of fields whose Galois group is a (possibly noncommutative) $p$-adic Lie group. This attempt to establish a non-commutative
version of the work of Iwasawa has led to a plethora of new problems and conjectures (for instance, see \cite{BH, CFKSV, FK, Ho, Ho2, OcV03, V02, V03Crelle}, where this list is far from being exhaustive). Indeed, Venjakob has previously predicted that an analogous asymptotic formula should hold for class groups over a false-Tate
extension (see \cite[p.\ 187]{V03Crelle}) which has recently been verified by Lei \cite{Lei} (also see \cite{LiangL}). Prior to this, Perbet \cite{Per} has also obtained results in this direction, although his methods do not allow
one to obtain a formula for the full $p$-class group but rather the $p^n$-torsion subgroup of the class group of the appropriate subfield $F_n$. We should also mention that even before the emergence of noncommutative Iwasawa theory, Cuoco-Monsky has studied the variation of class groups over a $\Zp^d$-extension (see \cite{CuoMo}).
The goal of this paper is to study the variation of even $K$-groups over $p$-adic Lie extensions which do not necessarily contain the cyclotomic $\Zp$-extension.


We now say a little on our approach. Before that, we quickly mention that in the class groups situation, one of the main tools is class field theory. In the works of Coates \cite{C72} and Ji-Qin \cite{JQ}, they have restricted themselves to the case of a cyclotomic $\Zp$-extension and the main reason of this restriction is so that they could rely on an Iwasawa theoretical version of Kummer theory (developed by Iwasawa \cite{Iw73}) to relate the even $K$-groups to the class groups situation. However, for our problem in hand, as we are concerned with $p$-adic Lie extensions which do not necessarily contain the cyclotomic $\Zp$-extension, we do not have this Kummer theoretical tool (and hence class field theory) at our disposal. Therefore, our approach takes a different route which we now explain. Let $i$ be a given integer $\geq 2$. Let $F_\infty$ be a $p$-adic Lie extension of $F$ which, for simplicity here, is assumed to be unramified outside $p$ (see the body of the paper for the general discussion). Building on the Quillen-Lichtenbaum Conjecture which is now a theorem thanks to the works of Rost-Voevodsky \cite{Vo, Wei09}, we can identify the Sylow $p$-subgroup of $K_{2i-2}(\Op_L)$ with the continuous cohomology group $H^2(G_{S_p}(L), \Zp(i))$ for every finite extension $L$ of $F$ contained in $F_\infty$. One is therefore reduced to studying the variation of these cohomology groups in the said $p$-adic Lie extension $F_\infty$. For this, we borrow an inspiration from Iwasawa's original works (\cite{Iw59, Iw73}), where he considered the inverse limit of the class groups with transition maps given by the norm maps on the class groups. In our situation,
we shall mimic this procedure by considering the inverse limit of the cohomology groups $H^2(G_{S_p}(L), \Zp(i))$ with transition maps given by the corestriction maps on the cohomology. The resulting limit module is in turn denoted by $H^2_{\Iw,S_p}(F_\infty/F, \Zp(i))$ which is called the (second) Iwasawa cohomology group. Therefore, the problem in hand is to analyse this Iwasawa cohomology group, which is now a module over the Iwasawa algebra of the Galois group of the $p$-adic Lie extension in question. As a start, we shall show that this module is torsion (see Proposition \ref{torsion H2}). We note that our torsionness result is obtained for a general $p$-adic Lie extension which does not necessarily contain the cyclotomic $\Zp$-extension and whose Galois group is not necessarily solvable. Indeed, when the $p$-adic Lie extension either contains the cyclotomic $\Zp$-extension or has a solvable Galois group, there are known methods in proving the said torsionness property (see \cite{BH, OcV03}; also see Remark \ref{torsion H2 remark} for discussion on these). Here our method of proof differs from these previous methods (see the proof of Proposition \ref{torsion H2} for details).  The next step is then to understand the descent of this Iwasawa cohomology group to an intermediate extension $L$. For this, we shall call upon an Iwasawa-cohomological version of a descent spectral sequence of Tate (obtained in \cite{FK, LimSh, Ne}) which tells us that \[H^2_{\Iw,S_p}(F_\infty/F, \Zp(i))_{\Gal(F_\infty/L)}\cong H^2(G_{S_p}(L), \Zp(i))\]
(see Proposition \ref{gal descent}).

We then combine the above descent observation with the various algebraic asymptotic formulas to obtain the variation of the even $K$-groups. When the $p$-adic Lie extension $F_\infty$ is commutative, we shall make use of the asymptotic formula of Cuoco-Monsky to obtain growth formulas for $K_{2i-2}(\Op_{F_n})[p^\infty]$ for appropriate intermediate subextensions $F_n$ of $F_\infty$ (see Theorem \ref{Zpd asymptotic arith}). Specializing to the cyclotomic $\Zp$-extension, we will see that our result recovers those as established by Coates, Jin and Qin. More precisely, we show that the Iwasawa invariants occurring in our asymptotic formulas, which comes from the Iwasawa cohomology groups, coincides with the Iwasawa invariants in the asymptotic formulas of Coates, Jin and Qin, which comes from certain Galois groups (see Remark \ref{Zp asymptotic arith rem)}).

 For a noncommutative $p$-adic Lie extension though, due to a lack of a precise structural theorem for modules over a noncommutative Iwasawa algebra, we can only obtain a growth formula for $K_{2i-2}(\Op_{F_n})[p^n]$ (see Theorem \ref{uniform asymptotic arith}) by appealing to the results of Perbet \cite{Per}. We also mention that in the noncommutative situation, one also needs to take some extra care due to the possible presence of ramified primes outside $p$ in the said noncommutative $p$-adic Lie extension. In the special case that the extension is a $\Zp^{d-1}\rtimes\Zp$-extension, we do have a precise asymptotic formula under certain extra hypotheses (see Theorem \ref{uniform asymptotic arith3}). For this latter result, we rely on a strategy of Lei \cite{Lei} and a technical result of Liang-Lim \cite{LiangL}, where such ideas were previously applied to study the growth of class groups over a $\Zp^{d-1}\rtimes\Zp$-extension. We take the opportunity here to formalize these ideas in the form of Proposition \ref{main alg theorem3}, as we believe that this result will surely have applications to Iwasawa
theory in other settings.

We end the introductional section by giving an outline of the paper.
In Section \ref{Iwasawa modules}, we review the notion of certain Iwasawa invariants. We also collect certain estimates on the intermediate coinvariants of a $\Zp\ps{G}$-module $M$ which will be required for our discussion on the variation of the even $K$-groups. Section \ref{Arithmetic preliminaries} is where we recall material from literature to describe the relation between the Sylow $p$-subgroups of the even $K$-groups and various \'etale/Galois cohomology groups. This paves the way for us to reduce the problem to a cohomological problem, which leads us to introduce the Iwasawa cohomology group in formal in Section \ref{Main results section}. It is here where we examine the module structure of the Iwasawa cohomology group. This is then combined with the discussion in the previous sections to yield our growth formulas in various classes of $p$-adic Lie extensions. Finally, in Section \ref{examples and remark}, we mention briefly some examples and classes of extensions that our results can be applied to. We also give some classes of extensions, where the second Iwasawa cohomology groups have positive $\mu$-invariants.

\subsection*{Acknowledgement}
The author would like to thank John Coates and Antonio Lei for their interest and comments. This research is supported by the
National Natural Science Foundation of China under Grant No. 11550110172 and Grant No. 11771164.

\section{Iwasawa modules} \label{Iwasawa modules}

\subsection{Iwasawa invariants}

We begin recalling the definition of certain Iwasawa invariants. Let $G$ denote a compact $p$-adic Lie group with no $p$-torsion. The completed group algebra of $G$ over $\Zp$ is given by
 \[ \Zp\ps{G} = \plim_U \Zp[G/U], \]
where $U$ runs over the open normal subgroups of $G$ and the inverse
limit is taken with respect to the canonical projection maps. It is well known that $\Zp\ps{G}$ is
a Noetherian Auslander regular ring (cf. \cite[Theorem 3.26]{V02} or \cite[Theorem A.1]{LimFine}) with no zero divisors (cf.\
\cite{Neu}). In particular, since the ring $\Zp\ps{G}$ is Noetherian with no zero divisors, it admits a skew field $Q(G)$ which is flat
over $\Zp\ps{G}$ (see \cite[Chapters 6 and 10]{GW} or \cite[Chapter
4, \S 9 and \S 10]{Lam}). This in turn enables us to define the notion of $\Zp\ps{G}$-rank of a finitely generated $\Zp\ps{G}$-module $M$, which is given by
$$ \rank_{\Zp\ps{G}}(M)  = \dim_{Q(G)} (Q(G)\ot_{\Zp\ps{G}}M). $$
The module $M$ is then said to be a
torsion $\Zp\ps{G}$-module if $\rank_{\Zp\ps{G}} (M) = 0$.
One can show that $M$ is torsion over $\Zp\ps{G}$ if and only if $\Hom_{\Zp\ps{G}}(M,\Zp\ps{G})=0$ (for instance, see \cite[Lemma 4.2]{LimFine}). In the event that the torsion $\Zp\ps{G}$-module $M$ satisfies $\Ext^1_{\Zp\ps{G}}(M,\Zp\ps{G})=0$, we say that $M$ is a pseudo-null $\Zp\ps{G}$-module.

For a finitely generated $\Zp\ps{G}$-module $M$, denote by
$M[p^\infty]$ the $\Zp\ps{G}$-submodule of $M$ consisting of elements
of $M$ which are annihilated by some power of $p$. Howson \cite[Proposition
1.11]{Ho2}, and independently Venjakob \cite[Theorem 3.40]{V02}), showed that there is a
$\Zp\ps{G}$-homomorphism
\[ \varphi: M[p^\infty]\lra \bigoplus_{i=1}^s\Zp\ps{G}/p^{\alpha_i},\] whose
kernel and cokernel are pseudo-null $\Zp\ps{G}$-modules, and where
the integers $s$ and $\alpha_i$ are uniquely determined. The $\mu_G$-invariant of $M$ is then defined to be $\mu_G(M) = \displaystyle
\sum_{i=1}^s\alpha_i$. We record the following basic lemma relating the rank and $\mu_G$-invariant.

\bl \label{mu and rank}
Let $M$ be a finitely generated $\Zp\ps{G}$-module.
Suppose that there is a $\Zp\ps{G}$-homomorphism
\[ \varphi: M[p^\infty] \longrightarrow \bigoplus_{i=1}^s\Zp\ps{G}  /p^{\alpha_i},\] whose
kernel and cokernel are pseudo-null $\Zp\ps{G}$-modules.
Then
\[\mu_{G}(M/p^n) =
n\rank_{\Zp\ps{G}}(M) + \sum_{i=1}^s\min\{n,\alpha_i\} \quad
\mbox{for}~ n\geq 1.
\]
In particular, if $M$ is a torsion $\Zp\ps{G}$-module, then we have $\mu_G(M)>0$ if and only if $\mu_{G}(M/p)>0$.
\el

\bpf
The equality is proven in \cite[Lemma 2.4.1]{LimMHG}. The final assertion of the lemma is immediate from this.
\epf

\subsection{The case of commutative $G$}

In this subsection, $G$ is taken to be a $p$-adic Lie group isomorphic to $\Zp^d$ for some $d\geq 1$. Write $G_n = G^{p^n}=(p^n\Zp)^d$. We shall collect certain results on the asymptotic behaviour of $M_{G_n}[p^\infty]$ for a $\Zp\ps{G}$-module $M$ that will be required for subsequent discussion. As a start, we recall the basic case $G=\Zp$.

\bp \label{Zp asymptotic}
 Suppose that $G=\Zp$ and that $M$ is a finitely generated $\Zp\ps{G}$-module. Assume that $M_{G_n}$ is finite for each $n$, where $G_n = p^n\Zp$. Then we have
 \[ \log_p\Big|M_{G_n}[p^\infty]\Big| = \mu_G(M)p^n + \lambda n + O(1)\]
 for a certain integer $\la$ independent of $n$.
\ep

\bpf
This is a special case of \cite[Proposition 5.3.17]{NSW}. Here $\lambda$ is precisely the $\lambda$-invariant of the module $M$ in the usual sense (see \cite[Definition 5.3.9]{NSW}) noting our hypothesis that $M_{G_n}$ is finite for each $n$.
\epf

The above algebraic result is fundamental in obtaining an asymptotic formula for the class groups in a $\Zp$-extension as was done by Iwasawa \cite{Iw59,Iw73}. Following this line of thought in their study of class groups in a $\Zp^d$-extension, Cuoco-Monsky established the following variant of the above formula.

\bp \label{Zpd asymptotic}
 Let $G=\Zp^d$ for some $d\geq 2$, and let $M$ be a finitely generated torsion $\Zp\ps{G}$-module. Suppose that $\rank_{\Zp}(M_{G_n}) = O(p^{(d-2)n})$. Then we have
 \[ \log_p\Big|M_{G_n}[p^\infty]\Big| = \mu_G(M)p^{dn} + l_0(M)np^{(d-1)n}+ O(p^{(d-1)n})\]
  for a certain integer $l_0(M)$ independent of $n$.
\ep

\bpf
This is \cite[Theorem 3.4]{CuoMo}. For the precise description of $l_0(M)$, we refer readers to \cite[Definition 1.2]{CuoMo}.
\epf

The above proposition suffices for our purposes when working over a multiple $\Zp^d$-extensions. However, when dealing with certain noncommutative $p$-adic Lie extensions (see Subsection \ref{ZprtimeZp}), we shall require the following variant which has no constraint on the module in hand but gives a less precise estimate.

\bp \label{Zpd asymptotic LiangLim}
 Let $G=\Zp^d$ and $M$ a finitely generated $($not necessarily torsion$)$ $\Zp\ps{G}$-module. Then we have
 \[ \log_p\Big|M_{G_n}[p^\infty]\Big| =\mu_{G}(M)p^{dn} + O(np^{(d-1)n}).\]
\ep

\bpf
See \cite[Theorem 2.4.1]{LiangL}.
\epf

\subsection{The case of $G$ being noncommutative} \label{noncom subsec}

 We come to the case of a noncommutative compact $p$-adic Lie group $G$. For simplicity, we shall assume that the group $G$ is uniform pro-$p$ in the sense of \cite[Definition 4.1]{DSMS}. Note that every compact $p$-adic Lie group contains an open normal subgroup which is a uniform pro-$p$ group by virtue of Lazard's theorem (see \cite[Corollary 8.34]{DSMS}). Therefore, one can always reduce consideration for a general compact $p$-adic Lie group to the case of a uniform pro-$p$ group which we will do throughout the paper.

Let $d$ denote the dimension of the uniform pro-$p$ group $G$. We write $G_n$ for the lower
$p$-series $P_{n+1}(G)$ which is defined recursively by $P_{1}(G) = G$, and
\[ P_{i+1}(G) = \overline{P_{i}(G)^{p}[P_{i}(G),G]}, ~\mbox{for}~ i\geq 1. \]
It follows from \cite[Theorem
3.6]{DSMS} that $G^{p^n} =
P_{n+1}(G)$ and that we have an equality $|G:P_2(G)| = |P_n(G):
P_{n+1}(G)|$ for every $n\geq 1$ (cf. \cite[Definition 4.1]{DSMS}). Hence we have  $|G:G_n| = |G:P_{n+1}(G)| = p^{dn}$.

With these notation in hand, we can now state the following result of Perbet.

\bp \label{Perbet thm}
 Let $G$ be a uniform pro-$p$ $p$-adic Lie group of dimension $d$ and $M$ a finitely generated $\Zp\ps{G}$-module. Then one has
  \[ \log_p\big|M_{G_n}/p\big| =\big(\rank_{\Zp\ps{G}}(M)+ \mu_G(M/p)\big)p^{dn} + O(p^{(d-1)n})\]
  and
 \[ \log_p\big|M_{G_n}/p^n\big| =\rank_{\Zp\ps{G}}(M)np^{dn} + \mu_G(M)p^{dn} + O(np^{(d-1)n}).\]
\ep

\bpf
See \cite[Th\'eor\`{e}me 2.1]{Per}.
\epf

A class of $p$-adic Lie groups $G$ that one is interested is the class of $G$ which contains a closed normal subgroup $H$ with the property that $G/H\cong\Zp$. As we shall see, in this context, the modules of interest usually are $\Zp\ps{G}$-modules which are finitely generated over $\Zp\ps{H}$-modules (see Proposition \ref{mu=0}). For such modules, they are clearly torsion over $\Zp\ps{G}$ and so have trivial $\Zp\ps{G}$-rank. Furthermore, it follows from \cite[Lemma 2.7]{Ho} that they have trivial $\mu_G$-invariant. Hence, in this context, Perbet's theorem reads $\log_p\big|M_{G_n}/p^n\big| = O(np^{(d-1)n})$. In view of this, we like to have more leverage on this $O(np^{(d-1)n})$ error term, which is the content of the next result.

\bp \label{main alg theorem2}
 Let $G$ be a uniform pro-$p$ $p$-adic Lie group of dimension $d$. Suppose that $G$ contains a closed normal subgroup $H$ with the property that $G/H\cong \Zp$. Let $M$ be a $\Zp\ps{G}$-module which is finitely generated over $\Zp\ps{H}$. Then we have
 \[ \log_p\left|M_{G_n}/p^n\right| \leq \rank_{\Zp\ps{H}}(M)np^{(d-1)n} + \mu_H(M)p^{(d-1)n} + O(np^{(d-2)n}).\]
\ep

\bpf
See \cite[Proposition 2.4]{LimUpper}.
\epf

\subsection{$G=\Zp^{d-1}\rtimes \Zp$} \label{ZprtimeZp}

In this subsection, $G$ is taken to be a uniform group which contains a closed normal subgroup $H\cong\Zp^{d-1}$ ($d\geq 2$) with $G/H\cong \Zp$. (We allow the case $G=\Zp^{d}$.) The goal of this subsection is to prove the following asymptotic behaviour of $M_{G_n}[p^\infty]$ which improves those in Subsection \ref{noncom subsec} under certain additional hypotheses on the module $M$.

\bp \label{main alg theorem3}
 Let $G$ be a compact pro-$p$ $p$-adic Lie group which contains a closed normal subgroup $H\cong \Zp^{d-1}$ such that $G/H\cong \Zp$. Let $M$ be a $\Zp\ps{G}$-module which is finitely generated over $\Zp\ps{H}$. Suppose further that $M_{G_n}$ is finite for every $n$. Then we have
 \[ \log_p\big|M_{G_n}[p^\infty]\big| = \rank_{\Zp\ps{H}}(M)np^{(d-1)n}  + O(p^{(d-1)n}).\]
\ep

The principle behind the above proposition has been fundamental in obtaining asymptotic class number formulas over a $\Zp^{d-1}\rtimes\Zp$-extension (see \cite{Lei, LiangL}). Here we formulate this precisely and give a proof of it in this form. The proof will take up the rest of this subsection.

We begin with some general remark.
Let $M$ be a finitely generated $\Zp\ps{G}$-module. For every subgroup $U$ of $G$, denote by $M(U)$ the $\Zp\ps{G}$-submodule of $M$ generated
by all elements of the form $(u-1)x$, where $u\in U$ and $x\in M$.
On the other hand, by restriction of scalar, we can view
$M$ as a $\Zp\ps{U}$-module. Denote by $I_UM$ the $\Zp\ps{U}$-submodule of $M$ generated by all elements of the form $(u-1)x$, where $u\in U$ and $x\in M$. Note that $I_UM$ needs not be a $\Zp\ps{G}$-module of $M$ if $U$ is not a normal subgroup of $G$.

\bl \label{coinvariant}
 For every subgroup $U$ of $G$, we have $I_UM\subseteq M(U)$. In the event that $U$ is a normal subgroup of $G$, we then have $I_UM = M(U)$ and
 \[M_U:= M/I_UM = M/M(U) \]
\el

\bpf
The inclusion is clear from the definition. Now suppose that $U$ is a normal subgroup of $G$. Then for each $g\in G, u\in U, x\in M$, we have
\[ g(u-1)x = (gug^{-1}-1)gx \in I_UM\]
by the normality of $U$. This establishes the lemma.
\epf

We return to the setting of this subsection.
Once and for all, fix a subgroup $\Ga$ of $G$ such that $\Ga$ maps isomorphically to $G/H$ under the natural quotient map $G\lra G/H$. In the event that $G=\Zp^d$, the group $\Ga$ is chosen such that $G= H\times \Ga$ (as an inner direct product). Either way, we have $G=H\Ga$ with $H\cap \Ga =1$. In other words, every element of $G$ can be expressed uniquely as $h\ga$ for $h\in H$ and $\ga\in \Ga$. For non-negative integers $m$ and $n$, we write $G_{m,n} = H_m\Ga_n$, where $H_m = H^{p^m}$ and $\Ga_n = \Ga^{p^n}$. In general, the $G_{m,n}$ are subgroups of $G$ but need not be normal subgroups. We also observe that $G_n = G^{p^n} = G_{n,n}$ and so, in particular, $G_n= G_{n,n}$ is a normal subgroup of $G$.

Clearly, $H_m$ is a characteristic subgroup of $H$ and so it follows that $H_m$ is a normal subgroup of $G$. By Lemma \ref{coinvariant}, we see that $M_{H_m} = M/M(H_m)$ has a natural $\Zp\ps{G/H_m}$-module structure. In particular, one can view $M_{H_m}$ as a $\Zp\ps{\Gamma}$-module, and consider
 \[\big(M_{H_m}\big)_{\Ga_n}  =M_{H_m} / I_{\Ga_n} M_{H_m}, \]
 where this a priori is an abelian group. On the other hand, as $H_m$ is a normal subgroup of $G$, we can define $\Zp\ps{G}$-module $\big(M/M(\Ga_n)\big)_{H_m}$.
 These various quotient-coinvariants are related as below.

\bl \label{quotient coinvariant}
For every integers $m,n$, we have
 \[\big(M_{H_m}\big)_{\Ga_n} = M/M(G_{m,n})
  =\big(M/M(\Ga_n)\big)_{H_m}. \]
  In particular, $\big(M_{H_m}\big)_{\Ga_n}$ can be endowed with a $\Zp\ps{G}$-module structure.
  \el

\bpf
Clearly, one has a natural surjection  $M_{H_m} \twoheadrightarrow M/M(G_{m,n})$
which factors through $\big(M_{H_m}\big)_{\Ga_n}$ to yield a surjection $f: \big(M_{H_m}\big)_{\Ga_n} \twoheadrightarrow M/M(G_{m,n})$. It remains to show that this is an injection. Suppose that $x\in M$ with $x\in M(G_{m,n})$. Then we have
\[ x =\sum_i h_i\ga_i(h_i'\ga_i'-1)y_i\]
for $h_i, h_i'\in H_m$, $\ga_i, \ga_i'\in \Ga_n$ and $y_i\in M$.
Rewriting, we have
\[x = \sum_i (h_i-1)\ga_i(h_i'\ga_i'-1)y_i + \sum_i \ga_i(h_i'-1)\ga_i'y_i \]
\[+ \sum_i (\ga_i-1)(\ga_i'-1)y_i +\sum_i (\ga_i'-1)y_i. \]
Note that the first two terms on the right lie in $M(H_m)$ and the final two terms on the right lie in $M(\Ga_n)$. This therefore shows that the class of $x$ in $\big(M_{H_m}\big)_{\Ga_n}$ is trivial which in turn implies that the map $f$ is injective. Hence we obtain the first identification of the lemma. The second can be shown similarly.
\epf

We record one more lemma which can be thought as an ``effective" version of Proposition \ref{Zp asymptotic}.

\bl \label{Gamma estimate}
Suppose that $\Ga\cong\Zp$.
Let $M$ be a finitely generated $\Zp\ps{\Ga}$-module with the properties that $M$ is finitely generated over $\Zp$ and that $M_{\Ga_n}$ is finite for every $n\geq 1$.
Let $f$ be a generator of the characteristic ideal of $M$. Let $n_0$ be an integer such that every irreducible distinguished polynomial that divides $f$ has degree $<p^{n_0-1}(p-1)$. Then
\[ \log_p\big|M_{\Ga_n}[p^\infty]\big| - \log_p\big|M_{\Ga_0}[p^\infty]\big| = \rank_{\Zp}(M)(n-n_0) + \log_p\big|M[p^\infty]_{\Ga_n}\big| -\log_p\big|M[p^\infty]_{\Ga_{n_0}}\big|.\]
for all $n\geq n_0$.   \el

\bpf
This is a special case of \cite[Proposition 4.6]{Lei} noting that $\mu_{\Ga}(M)=0$ as $M$ is finitely generated over $\Zp$.
\epf

We can now give the proof of Proposition \ref{main alg theorem3}.

\bpf[Proof of Proposition \ref{main alg theorem3}]
 Since $M$ is finitely generated over $\Zp\ps{H}$, it follows from \cite[Proposition 2.3]{CFKSV} that there exists a finite collection of $g_i\in \Zp\ps{G}$ such that each $\Zp\ps{G}/\Zp\ps{G}g_i$ is finitely generated over $\Zp\ps{H}$ and such that there is a surjection
 \[ \bigoplus_i \Zp\ps{G}/\Zp\ps{G}g_i \twoheadrightarrow M \]
 of $\Zp\ps{G}$-modules.

 As each $H_m$ is a normal subgroup of $G$, it follows from Lemma \ref{coinvariant} that $M_{H_m} = M/M(H_m)$ has a natural $\Zp\ps{G/H_m}$-module structure. In particular, one can view $M_{H_m}$ as a $\Zp\ps{\Gamma}$-module, and where the same can also be said for $\big(\Zp\ps{G}/\Zp\ps{G}g_i\big)_{H_m}$.
 By \cite[Corollary 5.4]{DL}, there exists an integer $t$ (independent of $m$) such that the $\Zp\ps{\Gamma}$-characteristic ideal of $\Big(\bigoplus_i \Zp\ps{G}/\Zp\ps{G}g_i\Big)_{H_m}$ factorizes into polynomials of degree $<t$. Since $M_{H_m}$ is a quotient of the direct sum of these modules, the same can also be said for the characteristic ideal of $M_{H_m}$. Now fix an integer $n_0$ such that $t< p^{n_0-1}(p-1)$. Also, for now, fix an arbitrary positive integer $m$. By  Lemma \ref{quotient coinvariant}, we have $\big(M_{H_m}\big)_{\Ga_n} = M/M(G_{m,n})$, where the latter is easily seen to be a quotient of $M_{G_{s}}$ for $s=\max\{m,n\}$. Since $M_{G_{s}}$ is finite by hypothesis, so is $\big(M_{H_m}\big)_{\Ga_n}$ for every $n$. Therefore, we may apply Lemma \ref{Gamma estimate} to $M_{H_m}$ to conclude that whenever $n\geq n_0$, one has
  \[ \begin{array}{l} \Big|\log_p\left|M/M(G_{m,n})\right| -  (n-n_0)\rank_{\Zp}(M_{H_m})\Big|\\
   \quad = \Big| \log_p\big|\big(M/M(\Ga_{n_0})\big)_{H_m}\big| + \log_p\big|M_{H_{m}}[p^{\infty}]_{\Ga_n}\big|  - \log_p\big|M_{H_{m}}[p^\infty]_{\Ga_{n_0}}\big| \Big| \\
   \quad \leq \log_p\left|\big(M/M(\Ga_{n_0})\big)_{H_m}\right| + 2 \log_p\big|M_{H_{m}}[p^\infty]\big|.
   \end{array} \]
   Since $m$ is arbitrary, the above inequality holds in particular when $m = n \geq n_0$. As a consequence, we have
   \[ \Big|\log_p\left|M_{G_n})\right|  - (n-n_0)\rank_{\Zp}(M_{H_n})\Big| \leq  \log_p\left|\big(M/M(\Ga_{n_0})\big)_{H_n}\right| + 2 \log_p\big|M_{H_{n}}[p^\infty]\big|.
   \]
   By Proposition \ref{Zpd asymptotic LiangLim}, the terms on the right hand side are $O(p^{(d-1)n})$ noting that $H$ has dimension $d-1$. On the other hand, it follows from  a result of Harris \cite[Theorem 1.10]{Har} that
   \[ \rank_{\Zp}(M_{H_n}) = \rank_{\Zp\ps{H}}(M)p^{(d-1)n} +O(p^{(d-2)n}), \]
   which in turn implies that
   \[ (n-n_0)\rank_{\Zp}(M_{H_n}) = \rank_{\Zp\ps{H}}(M)np^{(d-1)n} +O(p^{(d-1)n}). \]
   The conclusion of the proposition now follows by combining these estimates.
\epf

\section{Arithmetic preliminaries} \label{Arithmetic preliminaries}

In this section, we describe the relation between the algebraic $K$-groups and \'etale/Galois cohomology.

\subsection{\'Etale Cohomology and Galois cohomology}\label{etale Gal}

To begin with, we let $F$ be a number field, whose ring of integers is in turn denoted by $\Op_F$. Throughout the paper, $S$ will always denote a finite set of primes above $F$ which contains the set of primes above $p$ and the infinite primes. We then write $\Op_{F,S}$ for the ring of $S$-integers. We also write $S_\infty$ (resp., $S_p$) for the set of infinite primes (resp., primes above $p$) in $S$. Let $F_S$ be the maximal algebraic extension of $F$ unramified outside $S$ and write $G_S(F)$ for the Galois group $\Gal(F_S/F)$.
Denoting by $\mu_{p^n}$ the cyclic group generated by a primitive $p^n$-root of unity, we then write $\mu_{p^\infty}$ for the direct limit of the groups $\mu_{p^n}$. These have natural $G_S(F)$-module structures. Furthermore, for an integer $i\geq 2$, the $i$-fold tensor products $\mu_{p^n}^{\otimes i}$ and $\mu_{p^\infty}^{\otimes i}$ can be endowed with $G_S(F)$-module structure via the diagonal action. Therefore, we may speak of the Galois cohomology groups  $H^{k}\big(G_S(F), \mu_{p^n}^{\otimes i}\big)$ and  $H^{k}\big(G_S(F), \mu_{p^\infty}^{\otimes i}\big)$, noting that
\[H^{k}\big(G_S(F), \mu_{p^\infty}^{\otimes i}\big)\cong \ilim_n H^{k}\big(G_S(F), \mu_{p^n}^{\otimes i}\big). \]

On the other hand,  we can view $\mu_{p^n}^{\otimes i}$ as an \'etale sheaf over the scheme $\mathrm{Spec}(\Op_{F,S})$ in the sense of \cite[Chap.\ II]{Mi}, and consider the \'etale cohomology groups $H^{k}_{\acute{e}t}\big(\mathrm{Spec}(\Op_{F,S}), \mu_{p^n}^{\otimes i}\big)$.
By \cite[Chap.\ II, Proposition 2.9]{Mi}, the latter identifies with the Galois cohomology groups $H^{k}\big(G_S(F), \mu_{p^n}^{\otimes i}\big)$. Taking direct limit, we have
\[\ilim_n H^{k}_{\acute{e}t}\big(\mathrm{Spec}(\Op_{F,S}), \mu_{p^n}^{\otimes i}\big) \cong \ilim_n H^{k}\big(G_S(F), \mu_{p^n}^{\otimes i}\big)\cong H^{k}\big(G_S(F),\mu_{p^\infty}^{\otimes i}\big). \]
On the other hand, writing $\Zp(i) = \plim_n  \mu_{p^n}^{\otimes i}$ and taking inverse limit, we obtain
\[\plim_n H^{k}_{\acute{e}t}\big(\mathrm{Spec}(\Op_{F,S}), \mu_{p^n}^{\otimes i}\big) \cong \plim_n H^{k}\big(G_S(F), \mu_{p^\infty}^{\otimes i}\big) \cong H^{k}_{\mathrm{cts}}\big(G_S(F), \Zp(i)\big),\]
where $H^{k}_{\mathrm{cts}}(~,~)$ is the continuous cohomology group of Tate (see \cite[Chap 2, \S7]{NSW}), and where the second isomorphism is a consequence of \cite[Corollary 2.7.6 and Theorem 8.3.20(i)]{NSW}.

We end this subsection mentioning some convention on the notation that will be adhered for the remainder of the paper. For the \'etale cohomology groups, we shall always write $H^{k}_{\acute{e}t}\big(\Op_{F,S}, \Zp(i)\big)= \plim_n H^{k}_{\acute{e}t}\big(\mathrm{Spec}(\Op_{F,S}), \mu_{p^n}^{\otimes i}\big)$. While working with the continuous cohomology groups of Tate, we often drop ``cts" and write $H^{k}\big(G_S(F), \Zp(i)\big)= H^{k}_{\mathrm{cts}}\big(G_S(F), \Zp(i)\big)$. Finally, since the set of primes $S$ always contain $S_\infty$, in the event that $S=S_p\cup S_\infty$, we shall write $G_{S_p}(F)$ and $\Op_{F,S_p}$ for $G_{S_p\cup S_\infty}(F)$ and $\Op_{F,S_p\cup S_\infty}$ respectively.

\subsection{Algebraic $K$-theory} \label{Alg K-theory}

We now come to the $K$-theoretical aspects.
For a ring $R$ with identity, write $K_n(R)$ for the algebraic $K$-groups of $R$ in the sense of Quillen \cite{Qui73a, Qui73b} (also see \cite{Kol, Wei05, WeiKbook}). It is an elementary fact that $K_0(\Op_F)\cong\Z\oplus \mathrm{Cl}(F)$, where $\mathrm{Cl}(F)$ is the class group of $F$. By the theorem of Bass-Milnor-Serre \cite{BMS}, we have $K_1(\Op_F)\cong\Op_F^{\times}$, and so the structure of $K_1(\Op_F)$ can be read off from the Dirichlet's unit theorem. For the higher $K$-groups, we have the following deep result.

\bt \label{Borel}
The groups $K_{m}(\Op_F)$ are finitely generated for all $m\geq 0$. Furthermore, for $i\geq 2$, the groups $K_{2i-2}(\Op_F)$ are finite and
 \[ \rank_{\Z} K_{2i-1}(\Op_F) =\begin{cases} r_1(F)+r_2(F),  & \mbox{if $i$ is odd}, \\
r_2(F), & \mbox{if $i$ is even}. \end{cases} \]
Here $r_1(F)$ (resp., $r_2(F)$) is the number of real embeddings (resp., number of pairs of complex embeddings) of $F$.
\et

 \bpf
  Quillen \cite{Qui73b} showed that the $K$-groups are finitely generated. The asserted ranks of the $K$-groups are consequences of the calculations of Borel \cite{Bo}. We also note that the finiteness of $K_2(\Op)$ has previously been established by Garland \cite{Gar}.
 \epf

In \cite{Sou}, Soul\'e connected the higher $K$-groups with \'etale cohomology groups via the $p$-adic Chern class maps
\[ \mathrm{ch}_{i,k}^{(p)}: K_{2i-k}(\Op_F)\ot \Zp \lra H^k_{\acute{e}t}\left(\Op_{F,S_p}, \Zp(i)\right)\]
for $i\geq 2$ and $k =1,2$. (For the precise definition of these maps, we refer readers to loc.\ cit.) The famed Quillen-Lichtenbaum Conjecture predicts that these maps are isomorphisms. Thanks to the gallant efforts of many, we now know that this prediction is true.

\bt
 The $p$-adic Chern class maps are isomorphisms for $i\geq 2$ and $k =1,2$.
 \et

\bpf
Soul\'e first proved that these maps are surjective (see \cite[Th\'eor\`{e}me 6(iii)]{Sou}; also see the work of Dwyer and Friedlander \cite[Theorem 8.7]{DF}). It is folklore (for instance, see \cite[Theorem 2.7]{Kol}) that the asserted bijectivity is a consequence of the so-called norm residue isomorphism theorem (previously also known as the Bloch-Kato(-Milnor) conjecture; see \cite{BK, Mil}. Not to be confused with the other Bloch-Kato conjecture which is also known as the Tamagawa number conjecture; see \cite{BK2}). Later, Levine \cite{Lev} and Merkurjev-Suslin \cite{MS} established this said norm residue isomorphism theorem for the case $2i-k=3$. The full norm residue isomorphism theorem was eventually settled by the deep work of Rost and Voevodsky \cite{Vo} with the aid of a patch from Weibel \cite{Wei09}. Consequently, the $p$-adic Chern class maps are isomorphisms as asserted.
\epf

\bc \label{K2 = H2}
For $i\geq 2$, we have
\[ K_{2i-2}(\Op_F)[p^\infty] \cong H^2\big(G_{S_p}(F), \Zp(i)\big).\]
\ec

\bpf
Since $K_{2i-2}(\Op_F)$ is finite, we have $K_{2i-2}(\Op_F)[p^\infty]\cong K_{2i-2}(\Op_F)\ot\Zp$. By the preceding theorem, the latter is isomorphic to $H^2_{\acute{e}t}\big(\Op_{F,S_p}, \Zp(i)\big)$ which identifies with the corresponding Galois cohomology group as noted in Subsection \ref{etale Gal}.
\epf

In view of the preceding corollary, the study of $K_{2i-2}(\Op_F)[p^\infty]$ is reduced to the study of the cohomology groups $H^2\left(G_{S_p}(F), \Zp(i)\right)$. Occasionally, we may need to work with $H^2\left(G_{S}(F), \Zp(i)\right)$ for a larger set $S$ of primes, and so one needs to know the difference between these cohomology groups, which is the content of the next lemma.

\bl \label{change of S}
 Suppose that $S \supseteq S_p\cup S_\infty$. Then for $i\geq 2$, we have the following short exact sequence
\[ 0 \lra H^2_{\acute{e}t}\left(\Op_{F,S_p}, \Zp(i)\right)  \lra H^2_{\acute{e}t}\left(\Op_{F,S}, \Zp(i)\right)  \lra \bigoplus_{v\in S-S_p} H^1\left(k_v, \Zp(i-1)\right) \lra 0, \]
where $k_v$ is the residue field of $F_v$.
\el

\bpf
By the localization sequence of Soul\'e  (cf \cite[Proposition 2.2]{Kol} or \cite[Section III.3]{Sou}), we have the following exact sequence
\[ \bigoplus_{v\in S-S_p}H^0\left(k_v, \mu_{p^n}^{\otimes (i-1)}\right) \lra H^2_{\acute{e}t}\left(\Op_{F,S_p}, \mu_{p^n}^{\otimes i}\right)  \lra H^2_{\acute{e}t}\left(\Op_{F,S}, \mu_{p^n}^{\otimes i}\right) \]\[ \lra \bigoplus_{v\in S-S_p} H^1\left(k_v, \mu_{p^n}^{\otimes (i-1)}\right) \lra H^3_{\acute{e}t}\left(\Op_{F,S_p}, \mu_{p^n}^{\otimes i}\right). \]
As $G_{S_p}(F)$ has $p$-cohomological dimension $2$ (cf. \cite[Proposition 10.11.3]{NSW}), the final term is zero.
On the other hand, since $k_v$ is a finite field and $i\geq 2$, one has
\[ \plim_n H^0\left(k_v, \mu_{p^n}^{\otimes (i-1)}\right) = 0.\]
Hence we obtain a short exact sequence
\[0 \lra H^2_{\acute{e}t}\left(\Op_{F,S_p}, \mu_{p^n}^{\otimes i}\right)  \lra H^2_{\acute{e}t}\left(\Op_{F,S}, \mu_{p^n}^{\otimes i}\right) \lra \bigoplus_{v\in S-S_p} H^1\left(k_v, \mu_{p^n}^{\otimes (i-1)}\right) \lra 0. \]
The conclusion of the lemma now follows upon taking inverse limit of this short exact sequence.
\epf

We end with the following remark on the size of $H^1\left(k_v, \Zp(i-1)\right)$ which occurs in the exact sequence of Lemma \ref{change of S}.

\br \label{finite field remark} Let $i\geq 2$ be fixed. By \cite[Lemme 5]{Sou}, the group $H^1\big(k_v, \mu_{p^n}^{\otimes (i-1)}\big)$ is cyclic of order $\gcd(p^n, |k_v|^{i-1}-1)$. Since $k_v$ is a finite field, the group $\Gal(\bar{k}_v/k_v)$ has $p$-cohomological dimension 1. Consequently, the map \[H^1\big(k_v, \mu_{p^n}^{\otimes (i-1)}\big)\lra H^1\big(k_v, \mu_{p^m}^{\otimes (i-1)}\big)\] is surjective for $n\geq m$. It then follows that the order of $H^1\big(k_v, \mu_{p^n}^{\otimes (i-1)}\big)$ stabilizes and the above mentioned surjection is an isomorphism for $n\geq m\gg0$.
Hence $H^1\big(k_v, \Zp(i-1)\big)$ is a finite cyclic group of order $p^{\mathrm{ord}_p(|k_v|^{i-1}-1)}$.
\er

\section{Main results} \label{Main results section}

\subsection{Iwasawa cohomology groups} \label{Iw coh subsec}

Retain the notation as before. In particular, $F$ will denote a number field, and $S$ is a finite set of primes of $F$ which contains all the primes above $p$ and the infinite primes. For every extension $\mathcal{L}$ of $F$ contained in $F_S$, we write $G_S(\mathcal{L}) = \Gal(F_S/\mathcal{L})$. For $k \geq 0$ and $i\geq 2$, one defines the Iwasawa cohomology groups
 \[H^k_{\Iw,S}\big(\mathcal{L}/F, \Zp(i)\big):= \plim_L H^k\big(G_S(L), \Zp(i)\big), \]
 where $L$ runs through all finite extension of $F$ contained in $\mathcal{L}$ and the transition maps are given by the corestriction maps on the cohomology. Note that if $\mathcal{L}/F$ is a finite extension, we have
 $H^k_{\Iw,S}\big(\mathcal{L}/F, \Zp(i)\big)= H^k\big(G_S(\mathcal{L}), \Zp(i)\big)$.

 A Galois extension $F_\infty$ of $F$ is said to be a pro-$p$ $p$-adic Lie extension of $F$ if its Galois group $G$ is pro-$p$ compact $p$-adic Lie group without $p$-torsion and it is unramified outside a finite set of primes. For a pro-$p$ $p$-adic Lie extension $F_\infty$ contained in $F_S$, the modules $H^k_{\Iw,S}\big(F_\infty/F, \Zp(i)\big)$ are finitely generated over $\Zp\ps{G}$ (cf. \cite[Proposition 4.1.3]{LimSh}). For the second Iwasawa cohomology group $H^2_{\Iw,S}(F_\infty/F, \Zp(i))$, we can say even more.

\bp \label{torsion H2}
Let $F_\infty$ be a pro-$p$ $p$-adic Lie extension of $F$ contained in $F_S$. Then for every $i\geq 2$, $H^2_{\Iw,S}(F_\infty/F, \Zp(i))$ is a torsion $\Zp\ps{G}$-module.
\ep

We emphasis that the proposition does not require $F_\infty$ to contain the cyclotomic $\Zp$-extension $F^\cyc$ of $F$. For the proof of the proposition, we require the following well-known fact (see \cite{Sch}).

\bl \label{H2=0} Suppose that $i\geq 2$. Let $L$ be a finite extension of $F$ contained in $F_S$. Then we have $H^2(G_S(L),\mu_{p^\infty}^{\ot i}) = 0$.
\el

\bpf
 For the convenience of the reader, we give an explanation of this here.
 By considering the long continuous cohomology sequence of
 \[ 0\lra \Zp(i) \lra \Qp(i) \lra \mu_{p^\infty}^{\ot i}\lra 0, \]
 one obtains the following exact sequence
  \begin{equation} \label{long exact seq}
  H^2\big(G_S(L),\Zp(i)\big) \lra H^2\big(G_S(L),\Qp(i)\big)\lra H^2\big(G_S(L),\mu_{p^\infty}^{\ot i}\big)
 \lra H^3\big(G_S(L),\Zp(i)\big).
   \end{equation}
 Therefore, the conclusion of the lemma will follow from this once we can show that $H^2(G_S(L),\Qp(i))$ and  $H^3\big(G_S(L),\Zp(i)\big)$ vanish. Indeed, as $G_S(L)$ has $p$-cohomological dimension 2 (cf. \cite[Proposition 10.11.3]{NSW}), this gives the vanishing of $H^3\big(G_S(L),\Zp(i)\big)$. To see that $H^2(G_S(L),\Qp(i))$ vanishes, we first recall that as seen in Subsection \ref{Alg K-theory},
 the group $H^2(G_S(L), \Zp(i))$ is finite in view of Theorem \ref{Borel} and Corollary \ref{K2 = H2}. (Note that this latter finiteness property does not require the full power of the Quillen-Lichtenbaum conjecture. Indeed, the said finiteness is a well-known consequence of Borel's finiteness result on even $K$-groups and Soule's result on the surjectivity of the $p$-adic Chern class maps.) Combining these observations with (\ref{long exact seq}), we see that $H^2\big(G_S(L),\Qp(i)\big)$ sits between two torsion $\Zp$-modules. But $H^2\big(G_S(L),\Qp(i)\big)$ is also a $\Qp$-vector space, and so is $p$-torsionfree. Hence we must have $H^2\big(G_S(L),\Qp(i)\big)=0$. This therefore completes the proof of the lemma.
\epf

We require one more lemma.

\bl \label{torsion H2 van} Let $F_{\infty}$ be a pro-$p$
$p$-adic Lie extension of $F$ contained in $F_S$. Then the following statements are
equivalent.
\begin{enumerate}
\item[$(a)$] $H^2_{\Iw, S}(F_{\infty}/F, \Zp(i))$ is a torsion $\Zp\ps{G}$-module.

\item[$(b)$] $H^2(G_S(F_{\infty}),\mu_{p^\infty}^{\ot i}) = 0$.
\end{enumerate} \el

\bpf
 This is a special case of \cite[Lemma 7.1]{LimFine}.
\epf

We are in position to give the proof of Proposition \ref{torsion H2}.

\bpf[Proof of Proposition \ref{torsion H2}]
 By virtue of Lemma \ref{H2=0}, we have
  $H^2(G_S(L), \mu_{p^{\infty}}^{\ot i}) =0$ for every finite extension $L$ of $F$ contained in $F_\infty$. Taking direct limit over these $L$, we obtain  $H^2(G_S(F_\infty), \mu_{p^{\infty}}^{\ot i}) =0$. The conclusion of the proposition is now a consequence of this and Lemma \ref{torsion H2 van}.
\epf

\br \label{torsion H2 remark}
 We shall give two further proofs of Proposition \ref{torsion H2} but with some extra hypothesis. Although these proofs do not give the most general statement as in Proposition \ref{torsion H2}, we believe that it would be of interest to record them down.

(1) \underline{Suppose that $F_\infty$ contains the cyclotomic $\Zp$-extension of $F$.}

 \bpf
 Clearly, if the torsionness conclusion holds for an extension of $L$ over $F$, it holds over $F$. Thus, replacing the base field if necessary, we may assume that $F$ contains a $p$-th primitive root of unity. Then by hypothesis (1), $F_\infty$ contains $\mu_{p^{\infty}}$. This in turn implies that $G_S(F_\infty)$ acts trivially on $\mu_{p^{\infty}}^{\ot i}$.
 Hence we have $H^2(G_S(F_\infty), \mu_{p^{\infty}}^{\ot i}) = H^2(G_S(F_\infty), \Qp/\Zp)$ as abelian groups. But the latter is known to vanish (cf.\ \cite[Theorem 6.1]{OcV03}). The required conclusion is a consequence of this and Lemma \ref{torsion H2 van}.
 \epf

 Note that the above proof applies for other twists $i<2$ (still assuming that $F_\infty$ contains the cyclotomic $\Zp$-extension).

(2) \underline{Suppose that the Galois group $G=\Gal(F_\infty/F)$ is solvable.}

\bpf By Proposition \ref{gal descent}, one has
\[ H^2_{\Iw,S}\big(F_\infty/F, \Zp(i)\big)_{G} \cong H^2\big(G_S(F), \Zp(i)\big).\]
 Since $H^2\big(G_S(F), \Zp(i)\big)$ is finite and the group $G$ is assumed to be solvable, we may apply the main theorem of Balister-Howson \cite[Theorem in pp.\ 229]{BH} to conclude that
  $H^2_{\Iw,S}\big(F_\infty/F, \Zp(i)\big)$ is torsion over $\Zp\ps{G}$. \epf
\er

We end the subsection mentioning a version of Tate's descent spectral sequence for Iwasawa cohomology, which will be an important tool for all our proofs. In the statement, $T$ is taken to be a finitely generated  $\Zp$-module with a continuous $G_S(F)$-action and $H_i(U,-)$ is the usual homology group (for instance, see \cite[Definition 2.6.8]{NSW}).

\bp \label{gal descent}
Let $U$ be a closed normal subgroup of $G=\Gal(F_\infty/F)$ and write $L_\infty$ for the fixed field of $U$. Then we have a homological spectral sequence
 $$ H_r\big(U, H^{-s}_{\Iw,S}(F_\infty/F, T)\big)\Longrightarrow H^{-r-s}_{\Iw,S}\big(L_\infty/F, T\big). $$
 In particular, we have an isomorphism
\[ H^2_{\Iw,S}\big(F_\infty/F, T\big)_U \cong H^2_{\Iw,S}\big(L_\infty/F, T\big).\]
\ep

\bpf
Had $F_\infty//F$ being a finite extension, this is essentially the Tate spectral sequence (for instance, see \cite[Theorem 2.5.3]{NSW}). Thankfully, in the general context of the proposition, we have \cite[Proposition 1.6.5]{FK} or \cite[Theorem 3.1.8]{LimSh} (also see \cite[Proposition 4.2.3]{Ne}) to call upon.
 The final isomorphism in the proposition follows from reading off the initial $(0,-2)$-term of the spectral sequence.
\epf

\br \label{gal descent remark}
For our discussion, we usually take $T=\Zp(i)$ for some $i\geq 2$. We also note that the spectral sequence, and hence the asserted isomorphism of the second Iwasawa cohomology, is valid if $U$ is a finite group which is not necessarily of $p$-power order (see \cite[Theorem 2.5.3]{NSW} or \cite[Theorem 3.1.8]{LimSh}). This fact will also be utilized in some of our later discussion.
\er

We now give an application of Proposition \ref{gal descent} (we will make use of this said proposition again for subsequent discussion).

\bc \label{K=0}
Let $i\geq 2$ be given.
Let $F_\infty$ be a pro-$p$ $p$-adic Lie extension of $F$ contained in $F_S$. Suppose that for every prime $v$ of $F$ not dividing $p$ and ramified in $F_\infty$, the quantity $|k_v|^{i-1}-1$ is coprime to $p$, where $k_v$ is the residue field of $F_v$. If $K_{2i-2}(\Op_F)[p]=0$, then $K_{2i-2}(\Op_L)[p]=0$ for every finite extension $L$ of $F$ contained in $F_\infty$.
\ec

\bpf
Let $S$ be the set of primes of $F$ consisting precisely of the primes above $p$, the infinite primes and the ramified primes of $F_\infty/F$.
 By Corollary \ref{K2 = H2} and Lemma \ref{change of S}, we have the following short exact sequence
\[ 0\lra K_{2i-2}(\Op_{F})[p^\infty] \lra H^2(G_S(F),\Zp(i))\lra \bigoplus_{v\in S-S_p} H^1(k_{v}, \Zp(i-1)) \lra 0.\]
Taking Remark \ref{finite field remark} into account, it then follows from the hypothesis that $H^1(k_{v}, \Zp(i-1))=0$ for $v\in S-S_p$. Consequently, we have $H^2(G_S(F),\Zp(i))=0$ from this and the hypothesis that $K_{2i-2}(\Op_F)[p]=0$. By virtue of Proposition \ref{gal descent}, this in turn implies that \[H^2_{\Iw,S}\big(F_\infty/F,\Zp(i)\big)_{\Gal(F_\infty/F)}=0.\]
We may now apply Nakayama lemma (cf. \cite[Theorem in pp.\ 226]{BH} or \cite[Lemma 5.2.18]{NSW}) to conclude that $H^2_{\Iw,S}\big(F_\infty/F,\Zp(i)\big)=0$. By another application of Proposition \ref{gal descent}, we see that
\[ H^2(G_S(L),\Zp(i))\cong H^2_{\Iw,S}\big(F_\infty/F,\Zp(i)\big)_{\Gal(F_\infty/L)} =0.\]
Since $K_{2i-2}(\Op_{L})[p^\infty]$ is contained in $H^2(G_S(L),\Zp(i))$ by
Corollary \ref{K2 = H2} and Lemma \ref{change of S}, this yields the required assertion of the corollary.
\epf

\br
The preceding corollary applies to $p$-adic Lie extensions that are unramified outside $p$, since then the ramification assumption is vacuous.
In particularly, this applies to a $\Zp^d$-extension, as such an extension is unramified outside $p$ (cf. \cite[Theorem 1]{Iw73}). Specializing to the case of $d=1$, we therefore recover \cite[Theorem 8]{C72} and \cite[Theorem 3.3(1)]{JQ}.
\er

\subsection{Growth in $\Zp^d$-extensions}
Recall that a $\Zp^d$-extension is unramified outside $p$ (cf.\ \cite[Theorem 1]{Iw73}). Therefore, in this subsection, we may take $S= S_p\cup S_\infty$ which we shall do. We begin considering the case of a $\Zp$-extension. For a cyclotomic $\Zp$-extension, this was first proved by Coates \cite[Theorem 9]{C72} for $K_2$, and subsequently by Ji-Qin \cite[Theorem 3.3(2)]{JQ} for the higher even $K$-groups. We now show that one has an analogue formula for an arbitrary $\Zp$-extension.

\bt \label{Zp asymptotic arith}
Let $i\geq 2$ be given.
 Let $F_\infty$ be a $\Zp$-extension of $F$ and let $F_n$ be the intermediate subfields of $F_\infty/F$ with $\Gal(F_n/F)\cong \Z/p^n\Z$.
 Then we have
 \[ \log_p\Big|K_{2i-2}(\Op_{F_n})[p^\infty]\Big| = \mu_G\Big(H^2_{\Iw, S_p}\big(F_\infty/F, \Zp(i)\big)\Big) p^{n} + \la^{(i)} n + O(1)\]
 for some $\la^{(i)}$ independent of $n$.
\et

\bpf
By Proposition \ref{gal descent}, we have
\[H^2_{\Iw,S_p}\big(F_\infty/F, \Zp(i)\big)_{G_n} \cong H^2\big(G_{S_p}(F_n), \Zp(i)\big),\]
where the latter is finite. The conclusion of the theorem therefore follows immediately from an application of Proposition \ref{Zp asymptotic}.
\epf

\br \label{Zp asymptotic arith rem)}
We now describe how the above result recovers those in \cite[Theorem 9]{C72} and \cite[Theorem 3.3(ii)]{JQ}, where they consider a cyclotomic $\Zp$-extension $F^\cyc$ of a field $F$ which contains a primitive $p$th root of unity. In this context, one has $F^\cyc = F(\mu_{p^\infty})$. Let $\mathcal{M}$ be the maximal abelian pro-$p$ extension of $F^\cyc$ which is unramified outside $p$. In \cite[Theorem 9]{C72} and \cite[Theorem 3.3(ii)]{JQ}, the Iwasawa invariants in their asymptotic formulas are given by the $\mu$- and $\la$-invariants of $\Gal(\mathcal{M}/F^\cyc)$. We shall show that these invariants coincide with those of our Iwasawa cohomology group $H^2_{\Iw,S_p}\big(F^{\cyc}/F, \Zp(i)\big)$. (Of course, the fact that the invariants coincide follows immediately from comparing the asymptotic formulas of our with those of Coates and Jin-Qin. But we thought it worthwhile to give the following more conceptual argument.)
As a start, observe that $\Gal(\mathcal{M}/F^\cyc)\cong H^1(G_{S_p}(F^\cyc),\Qp/\Zp)^{\vee}$, where $(-)^{\vee}$ is the Pontryagin dual. On the other hand, since $F^\cyc$ contains all the $p$-power roots of unity, we may apply \cite[Lemma 2.5.1(c)]{ShRes} to conclude that
 \[ H^2_{\Iw,S_p}\big(F^{\cyc}/F, \Zp(i)\big) \cong H^2_{\Iw,S_p}\big(F^{\cyc}/F, \Zp(1)\big) \ot \Zp(i-1),\]
 where $\Zp(1)$ is the Tate module of $\mu_{p^{\infty}}$. The above identification in turn tells us that the modules $H^2_{\Iw,S_p}\big(F^{\cyc}/F, \Zp(i)\big)$ and  $H^2_{\Iw,S_p}\big(F^{\cyc}/F, \Zp(1)\big)$ have the same $\mu_G$-invariant and $\la$-invariant. We shall now prove that the $\mu_G$-invariant and $\la$-invariant of $H^2_{\Iw,S_p}\big(F^{\cyc}/F, \Zp(1)\big)$ coincide with those of  $\Gal(\mathcal{M}/F^\cyc)$.

By considering the low degree terms of the spectral sequence
\[ \Ext^r_{\Zp\ps{G}}\Big(H^{3-s}_{\Iw,S_p}\big(F^\cyc/F, \Zp(1)\big),\Zp\ps{G}\Big)\Longrightarrow  H^{3-r-s}\big(G_{S_p}(F^\cyc),\Qp/\Zp\big)^{\vee}\]
(cf.\ \cite[1.6.12]{FK},  \cite[Theorem 4.5.1]{LimSh} or \cite[Theorem 5.4.5]{Ne}), we obtain an exact sequence
\[ 0 \lra \Ext^1_{\Zp\ps{G}}\Big(H^{2}_{\Iw,S_p}(F^\cyc/F, \Zp(1)),\Zp\ps{G}\Big)\lra  H^{1}(G_{S_p}(F^\cyc),\Qp/\Zp)^{\vee} \]\[\lra \Ext^0_{\Zp\ps{G}}\Big(H^{1}_{\Iw,S_p}(F^\cyc/F, \Zp(1)),\Zp\ps{G}\Big) \]
Since $\Gal(\mathcal{M}/F^\cyc)\cong H^1(G_{S_p}(F^\cyc),\Qp/\Zp)^{\vee}$ and the $\Ext^0$-term is reflexive by \cite[Corollary 5.1.3]{NSW}, we see that $\Ext^1_{\Zp\ps{G}}\Big(H^{2}_{\Iw,S_p}(F^\cyc/F, \Zp(1)),\Zp\ps{G}\Big)$ is isomorphic to the torsion submodule of $\Gal(\mathcal{M}/F^\cyc)$. By \cite[Proposition 5.5.13]{NSW}, one has an pseudo-isomorphism
\[ \Ext^1_{\Zp\ps{G}}\Big(H^{2}_{\Iw,S_p}(F^\cyc/F, \Zp(1)),\Zp\ps{G}\Big) \sim H^{2}_{\Iw,S_p}\big(F^\cyc/F, \Zp(1)\big)^{\iota},\]
where $M^\iota$ is the $\Zp\ps{G}$-module with underlying abelian group $M$ with $G$-action given by $g\cdot_{\iota} x = g^{-1}x$. It is straightforward to verify that $\mu_G$-invariant and $\la$-invariant are invariant under the $\iota$-action. Hence, in conclusion, the modules $H^2_{\Iw,S_p}\big(F^{\cyc}/F, \Zp(i)\big)$ and $\Gal(\mathcal{M}/F^\cyc)$ share the same $\mu_G$-invariant and $\la$-invariant. This therefore shows that Theorem \ref{Zp asymptotic arith} recovers the results of  Coates, Ji and Qin.
\er

The next result considers the case of a $\Zp^d$-extension.

\bt \label{Zpd asymptotic arith}
Let $i\geq 2$ be given.
 Let $F_\infty$ be a $\Zp^d$-extension of $F$ and let $F_n$ be the intermediate subfields of $F_\infty/F$ with $\Gal(F_n/F)\cong (\Z/p^n\Z)^{\oplus d}$.
 Then we have
 \[ \log_p\Big|K_{2i-2}(\Op_{F_n})[p^\infty]\Big| = \mu_G\big(H^2_{\Iw, S_p}\big(F_\infty/F, \Zp(i)\big)\big) p^{dn} + l_G^{(i)}np^{(d-1)n} + O(p^{(d-1)n}) \]
 for some integer $l_G^{(i)}$ independent of $n$.
\et

\bpf
Proposition \ref{torsion H2} tells us that $H^2_{\Iw,S_p}\big(F_\infty/F, \Zp(i)\big)$ is torsion over $\Zp\ps{G}$. On the other hand, Proposition \ref{gal descent} yields
\[H^2_{\Iw,S_p}\big(F_\infty/F, \Zp(i)\big)_{G_n} \cong H^2\big(G_{S_p}(F_n), \Zp(i)\big),\]
where the latter is finite. Therefore, the hypotheses in Proposition \ref{Zpd asymptotic} are satisfied, and so we may apply it to obtain the conclusion of the theorem.
\epf

\subsection{Growth in noncommutative extensions}

We now come to the situation of a noncommutative $p$-adic Lie extension. Let $F_\infty$ be a $p$-adic Lie extension of $F$ unramified outside a finite set of primes, whose Galois group $G$ is a uniform pro-$p$ group of dimension $d$. Throughout this section, $S$ is taken to be the set of primes of $F$ consisting precisely of the primes above $p$, the infinite primes and the ramified primes of $F_\infty/F$. In particular, if $v\in S-S_p$, then $v$ is ramified in $F_\infty/F$.

\bt \label{uniform asymptotic arith}
Let $i\geq 2$ be given.
  Suppose that $F_\infty$ is a $p$-adic Lie extension of $F$ unramified outside a finite set of primes, whose Galois group $G$ is a uniform pro-$p$ group of dimension $d$.  Writing $F_n$ for the intermediate subfields of $F_\infty/F$ with $\Gal(F_\infty/F_n) = G_n$, we have
 \[ \log_p\Big|K_{2i-2}(\Op_{F_n})[p^n]\Big| = \mu_G\big(H^2_{\Iw,S}\big(F_\infty/F, \Zp(i)\big)\big) p^{dn} +  O(np^{(d-1)n}). \]
\et

 For the proof of the theorem, we require two more lemmas. The first of which is an elementary number theoretical observation.

\bl
 Let $a$ and $b$ be positive integers. Suppose that $a = \ord_p(b-1)$. Then for every $n\geq 1$, we have $a+n = \ord_p(b^{p^n}-1)$.
\el

\bpf
For a lack of proper reference, we shall supply a proof here. The proof is by induction on $n$. Suppose that $n=1$. Then we have
\[ b^p-1 = (b-1)(1+ b +\cdots+ b^{p-1}).\]
Since $a = \ord_p(b-1)$ by hypothesis, it remains to show that
$\ord_p(1+ b +\cdots+ b^{p-1}) =1$. Note that
\[1+ b +\cdots+ b^{p-1} = p +(b-1)+\cdots +(b^{p-1}-1). \]
Therefore, if $a\geq 2$, then all the terms $b-1$,..., $b^{p-1}-1$
are divisible by $p^2$ and so we do have $\ord_p(1+ b +\cdots+ b^{p-1}) =1$ in this case. Now, suppose that $a=1$. Then, upon rewriting, we have
\[p+(b-1)+\cdots +(b^{p-1}-1) =p+ (b-1)\big(1+ (b+1) + \cdots +(b^{p-2}+\cdots +1)\big). \]
Since $1 = \ord_p(b-1)$ by hypothesis and our assumption, we are reduced to showing that
\[\ord_p\big(1+ (b+1) + \cdots +(b^{p-2}+\cdots +1)\big) \geq 1.\]
But since $b\equiv 1$ (mod $p$), we have
\[1+ (b+1) + \cdots +(b^{p-2}+\cdots +1) \equiv 1+ 2 +\cdots + (p-1) \equiv p(p-1)/2 \equiv 0 ~(\mathrm{mod}~p).\]
Hence this proves the lemma for $n=1$.

Suppose that $n\geq 1$ and that $a+n = \ord_p(b^{p^n}-1)$. Then observe that
\[b^{p^{n+1}}-1 = (b^{p^n}-1)(1+ b^{p^n} + \cdots +b^{p^n(p-1)}).\]
Thus, it suffices to show that
\begin{equation} \label{num} \ord_p\big((1+ b^{p^n} + \cdots +b^{p^n(p-1)})\big) = 1. \end{equation}
But by our induction hypothesis, we have $b^{p^n} \equiv 1~(\mathrm{mod}~p^{a+n}$), and so
\[ 1+ b^{p^n} + \cdots +b^{p^n(p-1)} \equiv p~(\mathrm{mod}~p^{a+n}). \]
Since $a+n\geq 2$, this gives (\ref{num}).
\epf

The next lemma is a consequence of the preceding one.

\bl \label{finite field bound}
Let $j\geq 1$ be fixed.
 Let $k$ be a finite field and $k_\infty$ a $\Zp$-extension of $k$. Writing $k_n$ the intermediate subfield of $k_\infty/k$ with $|k_n:k| =p^n$, we have
 \[ \mathrm{ord}_p(|k_n|^j -1) = O(n).\]
\el

\bpf
 If $\mathrm{ord}_p(|k_n|^j -1) = 0$ for every $n$, then the asserted estimate of the lemma is clearly valid. Suppose not, then without loss of generality, by base changing $k$, we may assume that $\mathrm{ord}_p(|k|^j -1) =a \geq 1$. For every $n$, we clearly have $|k_n|^j-1 = |k|^{jp^n}-1$. The conclusion now follows from this and the preceding lemma.
\epf

We can now give the proof of Theorem \ref{uniform asymptotic arith}.

\bpf[Proof of Theorem \ref{uniform asymptotic arith}]
By appealing to Proposition \ref{gal descent}, we see that
\[H^2_{\Iw,S}\big(F_\infty/F, \Zp(i)\big)_{G_n}/p^n \cong H^2\big(G_{S}(F_n), \Zp(i)\big)/p^n.\]
By Proposition \ref{Perbet thm}, and taking Proposition \ref{torsion H2} into account, we see that
\begin{equation} \label{H2 formula} \log_p\Big| H^2\big(G_{S}(F_n), \Zp(i)\big)/p^n \Big| = \mu_G\big(H^2_{\Iw,S}\big(F_\infty/F, \Zp(i)\big)\big) p^{dn} +  O(np^{(d-1)n}). \end{equation}

On the other hand, in view of Corollary \ref{K2 = H2} and Lemma \ref{change of S}, we have the short exact sequence
\[ 0\lra K_{2i-2}(\Op_{F_n})[p^\infty] \lra H^2(G_S(F_n),\Zp(i)) \lra \bigoplus_{w\in S(F_n)-S_p(F_n)} H^1(k_{n,w}, \Zp(i)) \lra 0,\]
where $k_{n,w}$ is the residue field of $F_{n,w}$. This in turn induces the following exact sequence
\[\bigoplus_{w\in S(F_n)-S_p(F_n)} H^1(k_w, \Zp(i-1))[p^n] \lra K_{2i-2}(\Op_{F_n})[p^\infty]\big/p^n \lra H^2(G_S(F_n),\Zp(i))\big/p^n\]
 \begin{equation} \label{K2 H2 exact sequence}
  \lra \bigoplus_{w\in S(F_n)-S_p(F_n)} H^1(k_{n,w}, \Zp(i-1))\big/p^n \lra 0.
 \end{equation}
Since $K_{2i-2}(\Op_{F_n})[p^\infty]$ is finite, one has
\begin{equation} \label{K2 finite formula} \log_p\Big|K_{2i-2}(\Op_{F_n})[p^n]\Big| = \log_p\Big| K_{2i-2}(\Op_{F_n})[p^\infty]\big/p^n\Big|.\end{equation}

The required estimate of the theorem will follow from (\ref{H2 formula}), (\ref{K2 H2 exact sequence}) and (\ref{K2 finite formula}) once we show that \begin{equation} \label{local bounds outside p} \log_p\left|\bigoplus_{w\in S(F_n)-S_p(F_n)} H^1(k_{n,w}, \Zp(i-1))\right| = O(np^{(d-1)n}).\end{equation}
Writing
\[ \bigoplus_{w\in S(F_n)-S_p(F_n)} H^1(k_{n,w}, \Zp(i-1)) = \bigoplus_{v\in S-S_p} \bigoplus_{w|v} H^1(k_{n,w}, \Zp(i-1)),\]
we are reduced to bounding $\bigoplus_{w|v} H^1(k_{n,w}, \Zp(i-1))$ for each $v\in S-S_p$. A combination of Remark \ref{finite field remark} and Lemma \ref{finite field bound} tells us that
\begin{equation} \label{local bound outside p} \log_p\big|H^1(k_{n,w}, \Zp(i-1))\big| = O(n). \end{equation}
On the other hand, since $F_\infty/F$ is ramified for primes $v\in S-S_p$, the decomposition group of $G$ at $v$ has dimension at least one. Thus, the number of primes of $F_n$ above each $v \in S-S_p$ is $O(p^{(d-1)n})$ (cf. \cite[Lemma 4.2]{HS}). Combining this observation with (\ref{local bound outside p}), we obtain the desired estimate (\ref{local bounds outside p}). This therefore concludes the proof of the Theorem.
\epf

We now consider a variant of Theorem \ref{uniform asymptotic arith}, where an asymptotic upper bound is obtained.

\bp \label{uniform asymptotic arith2}
Let $i\geq 2$ be given.
  Suppose that $F_\infty$ is a $p$-adic Lie extension of $F$ and $F_n$ is the intermediate subfield of $F_\infty/F$ with $\Gal(F_\infty/F_n) = G_n$. Assume further that the following statements are valid.

 $(a)$ $G$ contains a closed normal subgroup $H$ with $G/H\cong \Zp$.

 $(b)$ For each $v\in S-S_p$, the decomposition group of $G$ at $v$ has dimension 2.

 $(c)$ $H^2_{\Iw, S}\big(F_\infty/F, \Zp(i)\big)$ is finitely generated over $\Zp\ps{H}$.

 Then we have
 \[ \log_p\Big|K_{2i-2}(\Op_{F_n})[p^n]\Big| \leq \rank_{\Zp\ps{H}}\Big(H^2_{\Iw,S}\big(F_\infty/F, \Zp(i)\big)\Big) np^{(d-1)n} \quad\quad\quad\quad\quad\quad\]
 \[\quad\quad\quad \quad\quad\quad + ~\mu_{H}\Big(H^2_{\Iw,S}\big(F_\infty/F, \Zp(i)\big)\Big) p^{(d-1)n} +  O(np^{(d-2)n}). \]
\ep

\bpf
 The proof proceeds similarly to that in Theorem \ref{uniform asymptotic arith}, where we make use of Proposition \ref{main alg theorem2} in place of Proposition \ref{Perbet thm}. We also mention that condition (b) is used to show that the number of primes of $F_n$ above each $v \in S-S_p$ is $O(p^{(d-2)n})$ by a similar argument to that in \cite[Lemma 4.2]{HS}.
 This can then be combined with the estimate (\ref{local bound outside p}) to yield
 \[ \log_p\left|\bigoplus_{w\in S(F_n)-S_p(F_n)} H^1(k_{n,w}, \Zp(i-1))\right| = O(np^{(d-2)n}).\]
\epf

We end with the following case of a $\Zp^{d-1}\rtimes \Zp$-extension, where one can obtain a growth formula for $\log_p\Big|K_{2i-2}(\Op_{F_n})[p^\infty]\Big|$ rather than $\log_p\Big|K_{2i-2}(\Op_{F_n})[p^n]\Big|$ as in the previous two results.

\bt \label{uniform asymptotic arith3}
Let $i\geq 2$ be given.
  Suppose that $F_\infty$ is a $p$-adic Lie extension of $F$ and $F_n$ is the intermediate subfield of $F_\infty/F$ with $\Gal(F_\infty/F_n) = G_n$. Assume further that the following statements are valid.

 $(a)$ $G$ contains a closed normal subgroup $H$ with $H\cong \Zp^{d-1}$ and $G/H\cong \Zp$.

 $(b)$ For each $v\in S-S_p$, the decomposition group of $G$ at $v$ has dimension 2.

 $(c)$ $H^2_{\Iw, S}\big(F_\infty/F, \Zp(i)\big)$ is finitely generated over $\Zp\ps{H}$.

 Then we have
 \[ \log_p\Big|K_{2i-2}(\Op_{F_n})[p^\infty]\Big| = \rank_{\Zp\ps{H}}\Big(H^2_{\Iw,S}\big(F_\infty/F, \Zp(i)\big)\Big) np^{(d-1)n} +  O(p^{(d-1)n}). \]
\et

\bpf
 The proof proceeds similarly as before, where we now make use of Proposition \ref{Zpd asymptotic LiangLim}.
\epf

\section{Some further remarks} \label{examples and remark}

\subsection{On $p$-adic Lie extensions containing the cyclotomic $\Zp$-extension}

In this subsection, we consider a class of $p$-adic Lie extensions, namely those which contains the cyclotomic $\Zp$-extension, where Proposition \ref{uniform asymptotic arith2} applies. We first record the following.

\bp \label{mu=0}
Let $i \geq 2$ be given. Suppose that $F$ is a finite abelian extension of $\Q$. Let $F_\infty$ be a pro-$p$ $p$-adic Lie extension of $F$ which contains the cyclotomic $\Zp$-extension $F^\cyc$. Then $H^2_{\Iw,S}\big(F_\infty/F, \Zp(i)\big)$ is finitely generated over $\Zp\ps{H}$, where $H= \Gal(F_\infty/F^\cyc)$.
\ep

\bpf
 Observe that by Proposition \ref{gal descent}, one has
 \[ H^2_{\Iw,S}\big(F_\infty/F, \Zp(i)\big)_H \cong H^2_{\Iw,S}\big(F^\cyc/F, \Zp(i)\big).\]
 By Nakayama's lemma (for instance, see \cite[Lemma 3.2.1]{LimSh}), we are therefore reduced to showing that $H^2_{\Iw,S}\big(F^\cyc/F, \Zp(i)\big)$ is finitely generated over $\Zp$. As noted in Remark \ref{gal descent remark}, we have
 \[ H^2_{\Iw,S}\big(F^{\cyc}(\mu_p)/F, \Zp(i)\big)_{\Gal(F^{\cyc}(\mu_p)/F^\cyc)} \cong H^2_{\Iw,S}\big(F^\cyc/F, \Zp(i)\big).\]
Since $F$ is abelian over $\Q$, so is $F(\mu_p)$. Thus, it suffices to show the $\Zp$-finite generation property of $H^2_{\Iw,S}\big(F^\cyc/F, \Zp(i)\big)$ under the assumption that $F$ is a finite abelian extension of $\Q$ which contains $\mu_p$. In particular, under this said assumption, we have $F^\cyc = F(\mu_{p^\infty})$. By an application of \cite[Lemma 2.5.1(c)]{ShRes}, we see that
 \[ H^2_{\Iw,S}\big(F^{\cyc}/F, \Zp(i)\big) \cong H^2_{\Iw,S}\big(F^{\cyc}/F, \Zp(1)\big) \ot \Zp(i-1),\]
 where $\Zp(1)$ is the Tate module of $\mu_{p^{\infty}}$. We are therefore reduced to proving the $\Zp$-finite generation property of $H^2_{\Iw,S}\big(F^{\cyc}/F, \Zp(i)\big)$. It is well-known that this said finite generation property is equivalent to the Iwasawa $\mu$-conjecture (see \cite{Iw73}), which in turn is a theorem of Ferrero-Washington \cite{FW} under the assumption that $F/\Q$ is abelian. Hence the proof of the proposition is completed.
\epf

\br
As can be seen from the proof, if the Iwasawa $\mu$-conjecture is known in general, one can then remove the ``abelian" hypothesis in Proposition \ref{mu=0}.
\er

Proposition \ref{mu=0} thus provides possible examples of $p$-adic Lie extensions, where one can apply Proposition \ref{uniform asymptotic arith2} and Theorem \ref{uniform asymptotic arith3}.

For imaginary quadratic fields $F$, many examples of $K_2(\Op_F)$ (also coined as the tame kernel) have been computed by Browkin and Gangl \cite{BrowHer}.
For instance, from the table given there, we see that $K_2\big(\Q(\sqrt{-4683})\big) \cong \Z/2\times \Z/2 \times \Z/3 \times \Z/37$.
Hence if $p$ is an odd prime $\neq 3, 37$, it follows from Corollary \ref{K=0} that $K_2(\Op_L)[p^\infty] =0$ for every finite extension $L$ of $\Q(\sqrt{-4683})$ contained in the (unique) $\Zp^2$-extension of $\Q(\sqrt{-4683})$. For $p = 3$ or $37$, it follows from Proposition \ref{mu=0} and Theorem \ref{uniform asymptotic arith3} that we have
\[ \log_p\Big|K_{2}(\Op_{F_n})[p^\infty]\Big| = \rank_{\Zp\ps{H}}\Big(H^2_{\Iw,S_p}\big(F_\infty/F, \Zp(2)\big)\Big) np^{n} +  O(1), \] where the $F_n$'s are the intermediate subfields of the $\Zp^2$-extension of $\Q(\sqrt{-4683})$ and $H= \Gal(F_\infty/F^{\cyc})$.

\subsection{Positive $\mu_G$-invariants of Iwasawa cohomology groups}
Of course, there exist $p$-adic Lie extensions, where the second Iwasawa cohomology groups can have positive $\mu_G$-invariants. In this subsection, we present a class of examples of these. We begin with a preliminary lemma.

\bl \label{mod p H2}
Let $F$ be a number field which contains $\mu_p$. Then for each $i\geq 2$, we have
\[ \dim_{\mathbb{F}_p} \big(H^2(G_{S_p}(F), \mu_p^{\ot i})\big) = \dim_{\mathbb{F}_p} \big(\mathrm{Cl}_{S_p}(F)[p]\big) + |S_p| -1.  \]
\el

\bpf
 Since the number field $F$ contains $\mu_p$, we have $H^2(G_{S_p}(F), \mu_p^{\ot i})\cong H^2(G_{S_p}(F), \mu_p)\ot\mu_p^{\ot (i-1)}$ which in turn implies that
 \[\dim_{\mathbb{F}_p} \big(H^2(G_{S_p}(F), \mu_p^{\ot i})\big) = \dim_{\mathbb{F}_p} \big(H^2(G_{S_p}(F), \mu_p)\big). \]
 On the other hand, from the Poitou-Tate sequence, we have
 \[ 0\lra \mathrm{Cl}_{S_p}(F)/p \lra H^2(G_{S_p}(F), \mu_p) \lra \bigoplus_{S_p}\Z/p \lra \Z/p\lra 0.\]
 The required conclusion now follows from a combination of these observations and noting that $\dim_{\mathbb{F}_p} \big(\mathrm{Cl}_{S_p}(F)[p]\big) = \dim_{\mathbb{F}_p} \big(\mathrm{Cl}_{S_p}(F)/p\big)$ by the finiteness of $\mathrm{Cl}_{S_p}(F)$.
\epf

The following is the main result of this subsection.

\bp
Let $i\geq 2$ and let $G=\Zp^d$ for a positive integer $d$. Suppose that $p$ is a prime such that $p>2d+1$. Then there exist infinitely many pairs $(F, F_{\infty})$, where $F$ is a finite cyclic extension of $\Q(\mu_p)$ and $F_{\infty}$ a $\Zp^d$-extension of $F$ such that $\mu_G\Big(H^2_{\Iw,S_p}\big(F_\infty/F, \Zp(i)\big)\Big)>0$.
\ep

\bpf
By \cite[Theorem 5.2]{Cuo84}, there exists a cyclic extension $F^{(1)}$ of $\Q(\mu_p)$ with a $\Zp^d$-extension $F^{(1)}_{\infty}$ of $F^{(1)}$ such that $\mu_G\left(X_{F^{(1)}_{\infty}}\right)>0$, where $X_{F^{(1)}_{\infty}}$ is the Galois group of the maximal abelian unramified pro-$p$ extension of $F^{(1)}_\infty$. We claim that
 \begin{equation} \label{mu>0 claim1}
\mu_G\Big(H^2_{\Iw,S_p}(F_{1,\infty}/F, \Zp(i))\Big)>0.
\end{equation}
By Lemma \ref{mu and rank} and Proposition \ref{torsion H2}, it suffices to show that \begin{equation} \label{mu>0 claim2}
\mu_G\Big(H^2_{\Iw,S_p}(F_{1,\infty}/F, \Zp(i))\big/p\Big)>0.
\end{equation}
 We first consider the case when there exists a prime of $F$ above $p$ which splits completely in $F_\infty/F$. Then writing $F^{(1)}_n$ for the intermediate subfield of $F^{(1)}_\infty/F^{(1)}$ fixed by $(p^n\Zp)^d$, it follows from Lemma \ref{mod p H2} that
\[\dim_{\mathbb{F}_p} \big(H^2(G_{S_p}(F^{(1)}_n), \mu_p^{\ot i})\big) \geq S_p(F^{(1)}_n) \geq p^{dn}, \]
which in turn implies that
\[\mu_G\Big(H^2_{\Iw.S_p}(F_{1,\infty}/F, \Zp(i))\big)/p\Big)\geq 1>0. \]
by Proposition \ref{Perbet thm} and taking Proposition \ref{torsion H2} into account. Now, suppose that every prime of $F^{(1)}$ above $p$ does not split completely in $F_\infty/F$. Then the decomposition group of $\Gal(F^{(1)}_{\infty}/F^{(1)})$ at each prime $v$ of $F^{(1)}$ above $p$, which we denote by $G_v$, has dimension $\geq 1$. It follows from \cite[Theorem 5.4]{Jannsen89} that there is an exact sequence
\[ \bigoplus_{v\in S_p} \mathrm{Ind}^{G_v}_G (\Zp) \lra X_{F^{(1)}_{\infty}} \lra Y_{F^{(1)}_{\infty}} \lra 0,\]
where $Y_{F^{(1)}_{\infty}}$ is the Galois group of the maximal abelian unramified pro-$p$ extension of $F^{(1)}_\infty$ at which the primes above $S_p$ splits completely.
Since the group $G_v$ has dimension $\geq 1$ and $\mathrm{Ind}^{G_v}_G (\Zp)$ is finitely generated over $\Zp\ps{G_v}$, it follows from \cite[Lemma 2.7]{Ho} that $\bigoplus_{v\in S_p} \mathrm{Ind}^{G_v}_G (\Zp)$ is $\Zp\ps{G}$-torsion with trivial $\mu_G$-invariant. Consequently, one has
\begin{equation} \label{mu >0} \mu_G\big(Y_{F^{(1)}_{\infty}}\big) = \mu_G\big(X_{F^{(1)}_{\infty}}\big)>0.
\end{equation}
On the other hand, from the Poitou-Tate sequence, we have
\[ 0\lra Y_{F^{(1)}_{\infty}}\Big/p \lra H^2_{\Iw,S_p}\left(F^{(1)}_{\infty}/F^{(1)}, \mu_p\right) \lra \bigoplus_{S_p}\mathrm{Ind}^{G_v}_G(\Z/p).\]
In view of Lemma \ref{mu and rank} and (\ref{mu >0}), this implies that
\[\mu_G\left(H^2_{\Iw}\big(F^{(1)}_{\infty}/F^{(1)}, \mu_p\big)\right) \geq \mu_G\Big(Y_{F^{(1)}_{\infty}}\big/p\Big)>0.\]
Since the number field $F^{(1)}$ contains $\mu_p$, we have \[H^2_{\Iw,S_p}(F^{(1)}_{\infty}/F^{(1)}, \mu_p^{\ot i})\cong H^2_{\Iw,S_p}(F^{(1)}_{\infty}/F^{(1)}, \mu_p)\ot\mu_p^{\ot (i-1)},\]
and so one also has
\[\mu_G\left(H^2_{\Iw,S_p}\big(F^{(1)}_{\infty}/F^{(1)}, \mu_p^{\ot i}\big)\right)>0.\]
By appealing to \cite[Theorem 5.2]{Cuo84} again, we can find a cyclic extension $F^{(2)}$ of $\Q(\mu_p)$ with a $\Zp^d$-extension $F^{(2)}_{\infty}$ of $F^{(2)}$ such that $\mu\left(X_{F^{(2)}_{\infty}}\right)>\mu\left(X_{F^{(1)}_{\infty}}\right)$. In particular, the pair $(F^{(1)}, F^{(1)}_\infty)$ and $(F^{(2)}, F^{(2)}_\infty)$ are distinct. By a similar argument as above, we see that $\mu_G\Big(H^2_{\Iw,S_p}(F^{(2)}_\infty/F^{(2)}, \Zp(i))\big)\Big)>0$. Continuing
the above process iteratively, we obtain infinite pairs $(F^{(m)}, F^{(m)}_\infty)$ satisfying the conclusion of the proposition.
\epf

We end the paper remarking that the arguments used in the paper can also be applied to obtain similar results for $K_{2i-2}(O_{F,S})$ for more general set $S$ of primes, noting Lemma \ref{change of S}.

\end{document}